\documentclass[a4paper,10pt,reqno, english]{amsart}  
\usepackage[utf8]{inputenc}
\usepackage[T1]{fontenc}
\usepackage{amsmath,amsthm}
\usepackage{amsfonts,amssymb,enumerate}
\usepackage{url,paralist}
\usepackage{mathtools}  
\usepackage[colorlinks=true,urlcolor=blue,linkcolor=red,citecolor=magenta]{hyperref}
\usepackage{enumerate}
\usepackage{anysize}
\usepackage[arrow,curve,matrix,tips,2cell]{xy}  
  \SelectTips{eu}{10} \UseTips
  \UseAllTwocells

\theoremstyle{plain}
\newtheorem{theorem}{Theorem}[section]
\newtheorem{lemma}[theorem]{Lemma}

\newtheorem{conjecture}[theorem]{Conjecture}
\newtheorem*{theorem*}{Theorem}

\newtheorem*{claim*}{Claim}

\theoremstyle{definition}

\newcommand{\R}{\mathbb{R}}
\newcommand{\Z}{\mathbb{Z}}
\newcommand{\F}{\mathbb{F}}
\newcommand{\E}{\mathrm{E}}
\newcommand{\B}{\mathrm{B}}
\newcommand{\pt}{\mathrm{pt}}
\newcommand{\dist}{\mathrm{dist}}
\newcommand\Sym{\mathfrak S}
\newcommand\sk{{\rm sk}}

\newcommand{\im}{\operatorname{im}}

\newcommand{\ind}{\operatorname{Index}}

\newcommand{\conve}{\operatorname{conv}}

\newcommand{\vertices}{\operatorname{vert}}
\newcommand{\relint}{\operatorname{relint}}

\begin{document}

\title[The Topological Transversal Tverberg theorem plus Constraints]
{The Topological Transversal Tverberg theorem plus Constraints}

\author[Blagojevi\'c]{Pavle V. M. Blagojevi\'{c}} 
\thanks{The research by Pavle V. M. Blagojevi\'{c} leading to these results has
        received funding from DFG via Berlin Mathematical School.
        Also supported by the grant ON 174008 of the Serbian Ministry of Education and Science.}
\address{Inst. Math., FU Berlin, Arnimallee 2, 14195 Berlin, Germany\hfill\break%
\mbox{\hspace{4mm}}Mat. Institut SANU, Knez Mihailova 36, 11001 Beograd, Serbia}
\email{blagojevic@math.fu-berlin.de} 
\author[Dimitrijevi\'c Blagojevi\'c]{Aleksandra S. Dimitrijevi\'c Blagojevi\'c}
\thanks{The research by Aleksandra Dimitrijevi\'c Blagojevi\'c leading to these results has
        received funding from the  grant ON 174008 of the Serbian Ministry of Education and Science.}
\address{Mat. Institut SANU, Knez Mihailova 36, 11001 Beograd, Serbia}
\email{aleksandra1973@gmail.com}
\author[Ziegler]{G\"unter M. Ziegler} 
\thanks{The research by G\"unter M. Ziegler received funding from DFG via the Research Training Group ``Methods for Discrete Structures'' and the Collaborative Research Center TRR~109 ``Discretization in Geometry and Dynamics.''}  
\address{Inst. Math., FU Berlin, Arnimallee 2, 14195 Berlin, Germany} 
\email{ziegler@math.fu-berlin.de}



\date{March 31, 2016; revised October 10, 2016\\\indent
To appear in ``Discrete and Intuitive Geometry -- László Fejes Tóth 100 Festschrift''  
			(G. Ambrus, I. Bárány, K. J. Böröczky, G. Fejes Tóth, J. Pach, eds.),
			Bolyai Society Mathematical Studies series, Springer}


\begin{abstract}
In this paper we use the strength of the constraint method in combination with a generalized Borsuk--Ulam type theorem and a cohomological intersection lemma to show how one can obtain many new topological transversal theorems of Tverberg type.
In particular, we derive a topological generalized transversal Van Kampen--Flores theorem and a topological transversal weak colored Tverberg theorem.
\end{abstract}

 \maketitle
 
\section{Introduction}
\label{sec : Introduction}

At the Symposium on Combinatorics and Geometry in Stockholm 1989, Helge Tverberg formulated the following conjecture that in a special case coincides with his famous 1966 result \cite[Thm.\,1]{Tverberg1966}.

\begin{conjecture}[The transversal Tverberg conjecture]
	\label{conjecture : transversal tverberg conjecture}
	Let 
	\begin{compactitem}[\ \ $\bullet$ ]
	\item $m$ and $d$ be integers with $0\leq m\leq d-1$,
	\item $r_0,\ldots,r_m\geq 1$ be integers, and 
	\item $N_0=(r_0-1)(d+1-m),\dots, N_m=(r_m-1)(d+1-m)$.
	\end{compactitem}	
Then for every collection of sets $X_0,\dots,X_m\subset\R^d$ with $|X_0|=N_0+1,\ldots,|X_m|=N_m+1$, there exist an $m$-dimensional affine subspace $L$ of $\R^d$
and  $r_{\ell}$ pairwise disjoint subsets $X_{\ell}^{1},\ldots,X_{\ell}^{r_{\ell}}$ of $X_{\ell}$, for $0\leq\ell\leq m$,  such that
	\[
		\conve (X_{0}^{1}) \cap L\neq\emptyset,\ldots, \conve (X_{0}^{r_0}) \cap L\neq\emptyset,
		\ \ldots \ ,
		\conve (X_{m}^{1})\cap L\neq\emptyset,\ldots,\conve (X_{m}^{r_m}) \cap L\neq\emptyset.
	\]  
\end{conjecture}

\noindent
For $m=0$ this conjecture is Tverberg's well-known theorem. 
Tverberg and Vre\'cica published the full conjecture in 1993 \cite{Tverberg1993}.
They proved that it also holds for $m=d-1$ \cite[Prop.\,3]{Tverberg1993}. For $m=1$ and arbitrary $d$ they verified the conjecture only in the following three cases: $r_0=1$, $r_1=1$, and $r_0=r_1=2$ \cite[Prop.\,1]{Tverberg1993}.

\begin{figure}[ht]
\centering
\includegraphics[scale=1.3]{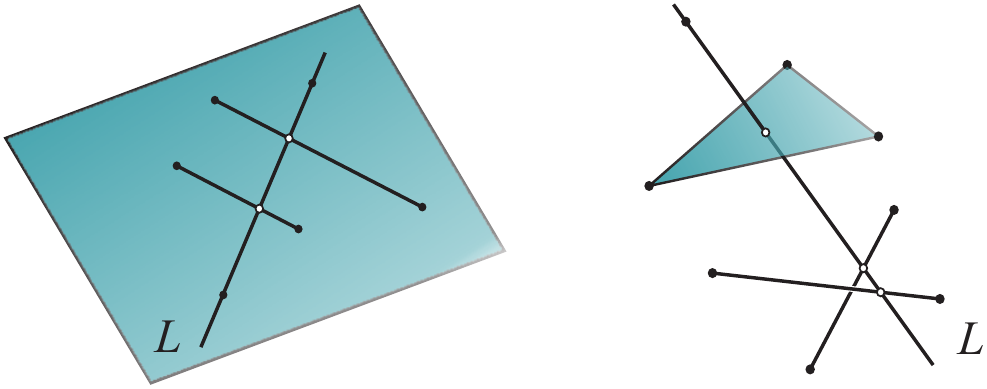}
\caption{\small The transversal Tverberg conjecture for $m=1$, $r_0=r_1=2$ and $d=2$ or $d=3$.}
\end{figure}

The classical Tverberg theorem from 1966 was extended to a topological setting by B\'ar\'any, Shlosman, and Sz\H{u}cs   \cite{Barany1981} in 1981. 
Similarly, it is natural to consider the following extension of the transversal Tverberg conjecture.

\begin{conjecture}[The topological transversal Tverberg conjecture]
	\label{conjecture : topological transversal tverberg conjecture}
	Let 
	\begin{compactitem}[\ \ $\bullet$ ]
	\item $m$ and $d$ be integers with $0\leq m\leq d-1$,
	\item $r_0,\ldots,r_m\geq 1$ be integers, and 
	\item $N_0=(r_0-1)(d+1-m),\dots, N_m=(r_m-1)(d+1-m)$.
	\end{compactitem}	
Then for every collection of continuous maps $f_0\colon\Delta_{N_0}\longrightarrow\R^d,\ldots, f_m\colon\Delta_{N_m}\longrightarrow\R^d$ there exist an $m$-dimensional affine subspace $L$ of $\R^d$ and  $r_{\ell}$ pairwise disjoint faces $\sigma_1^{\ell},\ldots,\sigma_{r_{\ell}}^{\ell}$  such that  
	\[
		f_{0}(\sigma_1^{0})\cap L\neq\emptyset, \ldots , f_{0}(\sigma_{r_{0}}^{0})\cap L\neq\emptyset,
		\quad \ldots \quad ,
		f_{m}(\sigma_1^{m})\cap L\neq\emptyset,\ldots , f_{m}(\sigma_{r_{m}}^{m})\cap L\neq\emptyset.
	\]
\end{conjecture}

\noindent
In 1999 using advanced methods of algebraic topology \v{Z}ivaljevi\'c \cite[Thm.\,4.8]{Zivaljevic1999} proved this conjecture for $d$ and $m$ odd integers and $r_0=\dots=r_m$ being an odd prime.
The topological transversal Tverberg conjecture was settled for $r_0=\dots=r_m=2$ by Vre\'cica \cite[Thm.\,2.2]{Vrecica2003} in 2003.
In 2007 Karasev \cite[Thm.\,1]{Karasev2007} established the topological transversal Tverberg conjecture in the cases when integers $r_0,\ldots,r_m$ are, not necessarily equal, powers of the same prime $p$ and the product $p(d-m)$ is even. 

\medskip
In the same paper Karasev \cite{Karasev2007} proved a colored topological transversal Tverberg's theorem \cite[Thm.\,5]{Karasev2007}, which for $m=0$ coincides with the colored Tverberg theorem of \v{Z}ivaljevi\'c and Vre\'cica \cite[Thm.\,pp.1]{Zivaljevic1992} and colored Tverberg theorem of type B of \v{Z}ivaljevi\'c and Vre\'cica \cite[Thm.\,4]{Vrecica1994}.
In 2011 Blagojevi\'c, Matschke and Ziegler gave yet another colored topological transversal Tverberg theorem \cite[Thm.\,1.3]{Blagojevic2011-01} that in the case $m=0$ coincides with their optimal colored Tverberg theorem \cite[Thm.\,2.1]{Blagojevic2009}.

\medskip
The existence of counterexamples to the topological Tverberg conjecture for non-primepowers, obtained by Frick \cite{Frick2015} \cite{Blagojevic2015-02} based on the remarkable work of Mabillard and Wagner \cite{Mabillard2014} \cite{Mabillard2015}, in particular invalidates Conjecture \ref{conjecture : topological transversal tverberg conjecture} in the case when $m=0$ and $r_0$ is not a prime power. 

\medskip
\subsubsection*{Acknowledgements.}
We are grateful to Florian Frick and to the referee 
for very good observations and useful comments. 

\section{Statement of the main results}

In 2014 Blagojevi\'c, Frick and Ziegler \cite{Blagojevic2014} introduced the ``constraint method,'' by which the topological Tverberg theorem implies almost all its extensions, which had previously been obtained as substantial independent results, such as the ``Colored Tverberg Theorem'' of \v{Z}ivaljevi\'c and Vre\'cica \cite{Zivaljevic1992} and the ``Generalized Van Kampen--Flores Theorem'' of Sarkaria \cite{Sarkaria1991-1} and Volovikov \cite{Volovikov1996-2}.
Thus the constraint method reproduced basically all Tverberg type theorems obtained during more than three decades with a single elementary idea.
Moreover, the constraint method in combination with the work of Mabillard and Wagner on the ``$r$-fold Whitney trick'' \cite{Mabillard2014} \cite{Mabillard2015} yields counterexamples to the topological Tverberg theorem for non-prime powers, as demonstrated by Frick \cite{Frick2015} \cite{Blagojevic2015-02}.

\medskip
In this paper we use the constraint method in combination with a generalized Borsuk--Ulam type theorem and a cohomological intersection lemma to show how one can obtain many new topological transversal theorems of Tverberg type.
We prove in detail a new generalized transversal van Kampen--Flores theorem and a new topological transversal weak colored Tverberg theorem.

\begin{theorem}[The topological generalized transversal Van Kampen--Flores theorem]
	\label{th : transversal generalized Van Kamplen Flores}
	Let 
	\begin{compactitem}[\ \ $\bullet$ ]
	\item $m$ and $d$ be integers with $0\leq m\leq d-1$,
	\item $r_0=p^{e_0},\ldots,r_m=p^{e_m}$ be powers of the prime $p$, where $e_0,\ldots,e_m\geq 0$ are integers,
	\item $N_0=(r_0-1)(d+2-m),\ldots, N_m=(r_m-1)(d+2-m)$,
	\item $k_0=\lceil \tfrac{r_0-1}{r_0}d \rceil ,\ldots, k_m=\lceil \tfrac{r_m-1}{r_m}d \rceil$, and
	\item $p(d-m)$ be even, or $m=0$.	
	\end{compactitem}
Then for every collection of continuous maps $f_0\colon\Delta_{N_0}\longrightarrow\R^d,\ldots, f_m\colon\Delta_{N_m}\longrightarrow\R^d$ there exist an $m$-dimensional affine subspace $L$ in $\R^d$ and  $r_{\ell}$ pairwise disjoint faces $\sigma_1^{\ell},\ldots,\sigma_{r_{\ell}}^{\ell}$ in the $k_{\ell}$-skeleton of $\Delta_{N_{\ell}}$, for $0\leq\ell\leq m$, such that  
	\[
		f_{0}(\sigma_1^{0})\cap L\neq\emptyset, \ldots , f_{0}(\sigma_{r_{0}}^{0})\cap L\neq\emptyset,
		\quad \ldots \quad ,
		f_{m}(\sigma_1^{m})\cap L\neq\emptyset,\ldots , f_{m}(\sigma_{r_{m}}^{m})\cap L\neq\emptyset.
	\]
\end{theorem}

\noindent
The special case $r_0=\cdots=r_m=2$ of the previous theorem is due to Karasev \cite[Cor.\,4]{Karasev2007}.

\medskip

In order to state the next result we recall the notion of a rainbow face. 
Suppose that the vertices of the simplex $\Delta_N$ are partitioned into color classes $\vertices(\Delta_N)=C_0\sqcup\cdots\sqcup C_k$.
The subcomplex $R_N:=C_0*\cdots * C_k\subseteq\Delta_N$ is called the \emph{rainbow complex}, 
that is, the subcomplex of all faces that have at most one vertex of each color class $C_0,\ldots,C_k$. 
Faces of $R_N$ are called \emph{rainbow faces}.

\begin{theorem}[The topological transversal weak colored Tverberg theorem]
	\label{th : topological transversal weak colored tverberg theorem}
	Let 
	\begin{compactitem}[\ \ $\bullet$ ]
	\item $m$ and $d$ be integers with $0\leq m\leq d-1$,
	\item $r_0=p^{e_0},\ldots,r_m=p^{e_m}$ be powers of the prime $p$, where $e_0,\ldots,e_m\geq 0$ are integers,
	\item $N_0=(r_0-1)(2d+2-m),\dots, N_m=(r_m-1)(2d+2-m)$,
	\item the vertices of the simplex $\Delta_{N_{\ell}}$, for every $0\leq m\leq d$, be colored by $d+1$ colors, where each color class has cardinality at most $2r_{\ell}-1$,
	\item $p(d-m)$ be even, or $m=0$.	
	\end{compactitem}
Then for every collection of continuous maps $f_0\colon\Delta_{N_0}\longrightarrow\R^d,\ldots, f_m\colon\Delta_{N_m}\longrightarrow\R^d$ there exist an $m$-dimensional affine subspace $L$ of $\R^d$ and $r_{\ell}$ pairwise disjoint rainbow faces $\sigma_1^{\ell},\ldots,\sigma_{r_{\ell}}^{\ell}$ in $\Delta_{N_{\ell}}$, for $0\leq\ell\leq m$, such that  
	\[
		f_{0}(\sigma_1^{0})\cap L\neq\emptyset, \ldots , f_{0}(\sigma_{r_{0}}^{0})\cap L\neq\emptyset,
		\quad \ldots \quad ,
		f_{m}(\sigma_1^{m})\cap L\neq\emptyset,\ldots , f_{m}(\sigma_{r_{m}}^{m})\cap L\neq\emptyset.
	\]
\end{theorem}

\medskip 
The proofs of Theorems~\ref{th : transversal generalized Van Kamplen Flores} and \ref{th : topological transversal weak colored tverberg theorem} are almost identical; see Sections \ref{subsec : GVKF} and \ref{subsec : WCT}.
The only difference occurs in the definition of the bundles $\xi_{\ell}$ and $\tau_{\ell}$, in \eqref{def of xi - 01} and \eqref{def of xi - 02}, and the bundle maps $\Phi_{\ell}$ for $0\leq\ell\leq m$; see Sections \ref{subsubsec: definition of function GVKF} and \ref{subsubsec: definition of function WCT}.
Using the same proof technique as for these theorems and modifying the bundles $\xi_{\ell}$ and bundle maps $\Phi_{\ell}$ using recipes from \cite[Lem.\,4.2]{Blagojevic2014}, one can also derive, for example, a topological transversal colored Tverberg theorem of type B, a topological transversal Tverberg theorem with equal barycentric coordinates, or mixtures of those.
The most general transversal Tverberg theorem that is produced by the constraint method can be formulated using the concept of “Tverberg unavoidable subcomplexes” \cite[Def.\,4.1]{Blagojevic2014}, as follows. 

Let $r \geq 2$, $d \ge 1$ and $N \geq r-1$ be integers, and let $f \colon \Delta_N \longrightarrow \R^d$ be a continuous map with at least one Tverberg $r$-partition, that is, a collection of $r$ pairwise disjoint faces $\sigma_1,\ldots,\sigma_r$ such that $f(\sigma_1)\cap\cdots\cap f(\sigma_r)\neq \emptyset$.
A subcomplex $\Sigma$ of the simplex $\Delta_N$ is \emph{Tverberg unavoidable with respect to $f$} 
if for every Tverberg partition $\{\sigma_1, \ldots, \sigma_r\}$ of $f$ there exists at least one face
$\sigma_i$ that lies in the subcomplex~$\Sigma$.%

\begin{theorem}[A constraint topological transversal Tverberg theorem]
	Let 
	\begin{compactitem}[\ \ $\bullet$ ]
	\item $m$ and $d$ be integers with $0\leq m\leq d-1$,
	\item $c_1,\ldots,c_m \geq 0$ be integer,
	\item $r_0=p^{e_0},\dots,r_m=p^{e_m}$ be powers of the prime $p$, where $e_0,\ldots,e_m\geq 0$ are integers,
	\item $N_0=(r_0-1)(d+1+c_1-m),\dots,N_t=(r_t-1)(d+1+c_m-m)$,
	\item $\Sigma_{i,j}$ be a Tverberg unavoidable subcomplex of the simplex $\Delta_{N_i}$ with respect to any continuous map $\Delta_{N_i}\longrightarrow\R^d$ for $1\leq i\leq m$ and $0\leq j\leq c_i$, assuming that $\Sigma_{i,0}=\Delta_{N_i}$, and
	\item $p(d-m)$ be even, or $m=0$.	
	\end{compactitem}
Then for every collection of continuous maps $f_0\colon\Delta_{N_0}\longrightarrow\R^d,\ldots, f_m\colon\Delta_{N_m}\longrightarrow\R^d$ there exist an $m$-dimensional affine subspace $L$ of $\R^d$ and $r_{\ell}$ pairwise disjoint faces $\sigma_1^{\ell},\ldots,\sigma_{r_{\ell}}^{\ell}$ that belong to the subcomplex  $\Sigma_{\ell,0}\cap\cdots\cap\Sigma_{\ell,c_i}$, for $0\leq\ell\leq m$, such that  
	\[
		f_{0}(\sigma_1^{0})\cap L\neq\emptyset, \ldots , f_{0}(\sigma_{r_{0}}^{0})\cap L\neq\emptyset,
		\quad \ldots \quad ,
		f_{m}(\sigma_1^{m})\cap L\neq\emptyset,\ldots , f_{m}(\sigma_{r_{m}}^{m})\cap L\neq\emptyset.
	\]
\end{theorem}

\section{A generalized Borsuk--Ulam type theorem and two lemmas}
\label{sec : methods}

In this section, we present the topological methods, developed in \cite{Blagojevic2011-01} and \cite{Karasev2007},  
that we will use in the proofs of Theorems \ref{th : transversal generalized Van Kamplen Flores} and \ref{th : topological transversal weak colored tverberg theorem}.
In particular, we will review and slightly modify a generalized Borsuk--Ulam type theorem \cite[Thm.\,4.1]{Blagojevic2011-01}, give an intersection lemma \cite[Lem.\,4.3]{Blagojevic2011-01} and recall the Euler class computation of Dol'nikov \cite[Lem.\,p.\,2]{Dolnikov1994}, \v{Z}ivaljevi\'c \cite[Prop.\,4.9]{Zivaljevic1999}, and Karasev \cite[Lem.\,8]{Karasev2007}.

\subsection{Fadell--Husseini index}

In 1988 Fadell and Husseini \cite{Fadell1988} introduced an ideal-valued index theory for the category of $G$-space, or more general for the category of $G$-equivariant maps to a fixed space with a trivial $G$-action.
We give an overview of the index theory adjusted to the needs of this paper.   

Let $G$ be a finite group, let $R$ be a commutative ring with unit, and let $B$ be a space with the trivial $G$ action.
For a $G$-equivariant map $p\colon X \longrightarrow B$ and a ring $R$, the \emph{Fadell--Husseini index} of $p$ is defined to be the kernel ideal of the map in the equivariant \v{C}ech cohomology with coefficients in the the ring $R$ induced by  $p$:
\begin{eqnarray*}
	\ind_{G}^{B}(p;R) & = & \ker\big(p^*\colon H^*(\E G\times_{G}B;R)\longrightarrow H^*(\E G\times _{G}X;R)\big)\\
	                  & = & \ker\big(p^*\colon H^*_G(B;R)\longrightarrow H^*_G(X;R)\big).
\end{eqnarray*}
The equivariant cohomology of a $G$-space $X$ is assumed to be the \v{C}ech cohomology of the Borel construction $\E G\times _{G}X$ associated to the space $X$.

\smallskip
\noindent
The basic properties of the index are:

\smallskip
\begin{compactitem}[\ \ $\circ$]

\item \textit{Monotonicity}: If $p\colon X \longrightarrow B$ and $q\colon Y \longrightarrow B$ are $G$-equivariant maps, and $f\colon X\longrightarrow Y$ is a $G$-equivariant map such that $p=q\circ f$, then
\[
\ind_{G}^{B}(p;R) \supseteq \ind_{G}^{B}(q;R).
\]

\item \textit{Additivity}: If $(X_1\cup X_2,X_1,X_2)$ is an excisive triple of $G$-spaces and $p\colon X_1\cup X_2 \longrightarrow B$ is a $G$-equivariant map, then
\[
\ind_{G}^{B}(p|_{X_1};R)\cdot\ind_{G}^{B}(p|_{X_2};R)\subseteq\ind_{G}^{B}(p;R).
\]
\item \emph{General Borsuk--Ulam--Bourgin--Yang theorem:} 
Let $p\colon X \longrightarrow B$ and $q\colon Y \longrightarrow B$ be $G$-equivariant maps, and let $f\colon X\longrightarrow Y$ be a $G$-equivariant map such that $p=q\circ f$.
If $Z\subseteq Y$ then
\[
\ind_{G}^{B}(p|_{f^{-1}(Z)};R)\cdot\ind_{G}^{B}(q|_{Y{\setminus}Z};R)\subseteq \ind_{G}^{B}(p;R).
\]
\end{compactitem}

\smallskip
\noindent
In the case when $B$ is a point and $p\colon X \longrightarrow B$ is a $G$-equivariant map we simplify notation and write $\ind_{G}^{B}(p;R)=\ind_{G}^{\pt}(X;R)$.
With this, the next property of the index can be formulated as follows.

\smallskip
\begin{compactitem}[\ \ $\circ$]
\item If $X$ is a $G$-space and $p\colon B\times X\longrightarrow B$ is the projection on the first factor, then
\[
\ind_{G}^{B}(p;R) = \ind_{G}^{\pt}(X;R)\otimes H^*(B;R).
\] 
\end{compactitem}

\subsection{A generalized Borsuk--Ulam type theorem}

The cohomology of the elementary abelian groups $(\Z/p)^e$, where $p$ is a prime and $e\geq 1$ is an integer, is given by 
\begin{equation*}
\begin{array}{llll}
H^*(\left( \Z/2\right) ^e;\F_2)        & =    & \F_2[t_1,\ldots,t_e],                                & \deg t_j=1                     \\
H^*(\left( \Z/p\right) ^e;\F_p)        & =    & \F_p[t_1,\ldots,t_e]\otimes \Lambda [u_1,\ldots,u_e], & \deg t_j=2,\deg u_i=1\text{ for $p$ odd.}
\end{array}
\end{equation*}
The following theorem and its proof is just a slight modification of \cite[Thm.\,4.1]{Blagojevic2011-01}.

\begin{theorem}[Borsuk--Ulam type theorem]
\label{th : BU}
Let
\begin{compactitem}[\ \ $\bullet$]
\item $G=(\Z/p)^e$ be an elementary abelian group where $p$ is a prime and $e\geq1$,
\item $B$ be a connected space with the trivial $G$-action,
\item $q\colon E \longrightarrow B$ be a $G$-equivariant vector bundle where all fibers carry the same $G$-representation,
\item $q|_{E^G}\colon E^G\longrightarrow B$ be the fixed-point subbundle of the vector bundle $q\colon E\longrightarrow B$,
\item $q|_{C}\colon C\longrightarrow B$ be its $G$-invariant orthogonal complement subbundle $(E=C\oplus E^G)$,
\item $F$ be the fiber of the vector bundle $q_{C}\colon C\longrightarrow B$ over the point $b\in B$,
\item $0\neq \alpha\in H^*(G;\F_p)$ be the Euler class of the vector bundle $F\longrightarrow \E G\times_G F\longrightarrow \B G$, and
\item $K$ be a $G$-CW-complex such that $\alpha\notin \ind_G^{\pt}(K;\F_p)$.
\end{compactitem} 
Assume that 
\begin{compactitem}[\ \ $\bullet$]
\item $\pi_1(B)$ acts trivially on $H^*(F;\F_p)$, and
\item we are given a $G$-equivariant map $\Phi \colon B\times K\longrightarrow E$ such that the following diagram commutes
\[
\xymatrix{
B\times K\ar[rr]^{\Phi} \ar[dr]_{q_1} && E=C\oplus E^G \ar[dl]^q\\
& B &
}
\]
where $q_1\colon B\times K\longrightarrow B$ is the projection on the first coordinate.
\end{compactitem} 

\smallskip
\noindent
Then for 
\[
S:=\Phi^{-1}(E^G)
\qquad\text{and}\qquad
T:=\Phi(S)=\im(\Phi)\cap E^G
\]
the following maps, induced by the projections $q_1$ and $q$, are injective:
\[
(q_1|_S)^*\colon H^*(B;\F_p) \longrightarrow  H^*_G(S;\F_p)
\qquad\text{and}\qquad
(q|_T)^*\colon H^*(B;\F_p)\longrightarrow  H^*(T;\F_p).
\]
\end{theorem}

\subsection{Two lemmas}
In this section we recall two facts: an intersection lemma from \cite[Lem.\,4.3]{Blagojevic2011-01} and the computation of a particular Euler class from \cite[Lem.\,8]{Karasev2007}.

\begin{lemma}[The intersection lemma]
\label{lemma : intersection lemma}
Let 
\begin{compactitem}[\ \ $\bullet$]
\item $k\geq 1$ be an integer, and $p$ a prime,
\item $B$ be an $\F_p$-orientable compact $m$-manifold,
\item $q\colon E\longrightarrow B$ be an $n$-dimensional real vector bundle whose mod-$p$ Euler class $e\in H^n(B;\F_p)$ satisfies $e^k\neq 0$, and 
\item $T_0,\ldots, T_k\subseteq E$ be compact subsets with the property that the induced maps 
	\[(q|_{T_i})^*\colon H^{m}(B;\F_p)\longrightarrow H^{m}(T_i;\F_p),
	\] 
	for all $0\leq i\leq k$, are injective.
\end{compactitem}
Then
\[
T_0\cap\cdots\cap T_k\neq \emptyset.
\]
\end{lemma}

\medskip
Let $G_n(\R^d)$ denote the Grassmann manifold of all $n$-dimensional subspaces in $\R^d$, and let $\gamma^n(\R^d)$ be the corresponding canonical vector bundle over $G_n(\R^d)$.
Furthermore, let $\widetilde{G}_n(\R^d)$ denote the oriented Grassmann manifold of all $n$-dimensional oriented subspaces in $\R^d$, and let $\widetilde{\gamma}^n(\R^d)$ be the corresponding canonical vector bundle over $\widetilde{G}_n(\R^d)$.
Then the ``forgetting orientation'' map  $\widetilde{G}_n(\R^d)\longrightarrow G_n(\R^d)$ is a double cover, and it induces a vector bundle map $\widetilde{\gamma}^n(\R^d)\longrightarrow\gamma^n(\R^d)$ that is an isomorphism on fibers.

\begin{lemma}[Euler classes of the canonical bundles of real Grassmannians]
\label{lemma : euler class}
Let $d$ and $m$ be positive integers with $0\leq m\leq d-1$, and let $p$ be a prime.
\begin{compactenum}[\ \ \rm (1)]
\item If $p=2$ and $\gamma:=\gamma^{d-m}(\R^d)$, then the $m$-th power of the Euler class of $\gamma$ does not vanish, that is
	\[
		0\neq e(\gamma)^m = w_{d-m}(\gamma)^m \in H^{(d-m)m}(G_{d-m}(\R^d);\F_2).
	\]
\item If $p$ is an odd prime, $d-m$ is even, and $\widetilde{\gamma}:=\widetilde{\gamma}^{d-m}(\R^d)$, then the $m$-th power of the mod-$p$ Euler class of $\widetilde{\gamma}$ does not vanish, that is
	\[
		0\neq e(\widetilde{\gamma})^m \in H^{(d-m)m}(\widetilde{G}_{d-m}(\R^d);\F_p).
	\]
\item If $p$ is an odd prime, $d-m$ is even, and $\gamma:=\gamma^{d-m}(\R^d)$, then the $m$-th power of the mod-$p$ Euler class of $\gamma$ does not vanish, that is
	\[
		0\neq e(\gamma)^m \in H^{(d-m)m}(G_{d-m}(\R^d);\F_p).
	\]
\end{compactenum}	
\end{lemma}

\noindent
The third part of the lemma is a consequence of the second part and the naturality property of the Euler class.
The case $p=2$ of this lemma was proved by Dol'nikov in \cite[Lem.\,p.\,112]{Dolnikov1994}.
For $p$ an odd prime, $d\geq 3$ an odd integer, and $d-m$ even the lemma was first proved  by \v{Z}ivaljevi\'c \cite[Prop.\,4.9]{Zivaljevic1999}.

\section{Proofs}

Now, combining the methods presented in Section \ref{sec : Introduction}, Theorem \ref{th : BU} and Lemmas \ref{lemma : intersection lemma} and \ref{lemma : euler class}, we prove our main results, Theorems \ref{th : transversal generalized Van Kamplen Flores} and \ref{th : topological transversal weak colored tverberg theorem}. 

\subsection{Proof of the topological generalized transversal Van Kampen--Flores theorem}
\label{subsec : GVKF}
Let $B:=G_{d-m}(\R^d)$ be the Grassmann manifold, and let $\gamma:=\gamma^{d-m}(\R^d)$ be the canonical bundle.
Without loss of generality we can assume that $m\geq 1$.
The proof of Theorem \ref{th : transversal generalized Van Kamplen Flores} is done in several steps.

\subsubsection{~~}

Fix an integer $0\leq\ell\leq m$, and define $K_{\ell}:=(\Delta_{N_{\ell}})^{\times r_{\ell}}_{\Delta(2)}$ to be the $r_{\ell}$-fold $2$-wise deleted product of the simplex $\Delta_{N_{\ell}}$. 
According to \cite[Lem.\,1]{Barany1981} the complex $K_{\ell}$ is an $(N_{\ell}-r_{\ell}+1)$-dimensional and $(N_{\ell}-r_{\ell})$-connected CW complex.
The symmetric group $\Sym_{r_{\ell}}$ acts freely on $K_{\ell}$ by permuting factors in the product, that is $\pi\cdot (x_1,\ldots,x_{r_{\ell}}):=	(x_{\pi^{-1}(1)},\ldots,x_{\pi^{-1}(r_{\ell})})$, for $\pi\in\Sym_{r_{\ell}}$ and $(x_1,\ldots,x_{r_{\ell}})\in K_{\ell}$.

Consider the regular embedding $\mathrm{reg}:(\Z/p)^{e_{\ell}}\longrightarrow \Sym_{r_{\ell}}$ of the elementary abelian group $(\Z/p)^{e_{\ell}}$, as explained in  \cite[Ex.\,2.7, p.\,100]{Adem2004}.
It is given by the left translation action of $(\Z/p)^{e_{\ell}}$ on itself. 
To every element $g\in (\Z/p)^{e_{\ell}}$ we associate the permutation $L_g\colon  (\Z/p)^{e_{\ell}}\longrightarrow  (\Z/p)^{e_{\ell}}$ from $\mathrm{Sym}((\Z/p)^{e_{\ell}})\cong\Sym_{r_{\ell}}$ given by $L_g(x)=g+x$.
Thus, the elementary abelian group $G_{\ell}:=(\Z/p)^{e_{\ell}}$ is identified with subgroup $\im(\mathrm{reg})$ of the symmetric group $\Sym_{r_{\ell}}$.
Consequently, $K_{\ell}$ is a free $G_{\ell}$-space.

Furthermore, let $\R^{r_{\ell}}$ be a vector space with the (left) action of the symmetric group $\Sym_{r_{\ell}}$ given by permutation of coordinates. 
Then the subspace $W_{r_{\ell}}:=\{(t_1,\ldots,t_{r_{\ell}})\in \R^{r_{\ell}} : \sum t_i=0\}$ is a $\Sym_{r_{\ell}}$-invariant subspace.
The group $G_{\ell}$ acts on both $\R^{r_{\ell}}$ and $W_{r_{\ell}}$ via the regular embedding.

Let $\tau_{\ell}$ be the the trivial vector bundle $B\times W_{r_{\ell}}\longrightarrow B$.
The action of $G_{\ell}$ on $W_{r_{\ell}}$ makes $\tau_{\ell}$ into a $G_{\ell}$-equivariant vector bundle.
Next,  $\gamma^{\oplus r_{\ell}}$ is also a $G_{\ell}$-equivariant vector bundle where the action is given by permutation of summands in the Whitney sum.
Then the vector bundle 
\begin{equation}
	\label{def of xi - 01}
	\xi_{\ell}:=\tau_{\ell}\oplus \gamma^{\oplus r_{\ell}}
\end{equation}
inherits the structure of a $G_{\ell}$-equivariant vector bundle via the diagonal action.
Let  $E(\cdot)$ denote the total space of a vector bundle.
Since the $G_{\ell}$ fixed point set of $W_{r_{\ell}}$ is just zero, that is $W_{r_{\ell}}^{G_{\ell}}=\{0\}$, the fixed point set of the total space of $\xi_{\ell}$ is
\[
 E(\xi_{\ell})^{G_{\ell}}=
 E(\tau_{\ell}\oplus \gamma^{\oplus r_{\ell}})^{G_{\ell}}\cong 
 E(\gamma^{\oplus r_{\ell}})^{G_{\ell}}\cong
 E(\gamma).
\]

\subsubsection{~~}
\label{subsubsec: definition of function GVKF}
We define a continuous $G_{\ell}$-equivariant bundle map $\Phi_{\ell}\colon B\times K_{\ell} \longrightarrow E(\xi_{\ell})$ as follows: For the point $(b,(x_1,\ldots,x_{r_{\ell}}))\in B\times K_{\ell}$ let
\begin{multline*}
\Phi_{\ell} (b,(x_1,\ldots,x_{r_{\ell}})):= \\ 
\qquad\qquad\qquad\Big(b,
\big(\dist(x_1, \sk_{k_{\ell}}(\Delta_{N_{\ell}}))-a(x_1,\ldots,x_{r_{\ell}}),\ldots, 
	 \dist(x_{r_{\ell}}, \sk_{k_{\ell}}(\Delta_{N_{\ell}}))-a(x_1,\ldots,x_{r_{\ell}}) \big)\Big) \oplus \\
\big((q_{b}\circ f_{\ell})(x_1) \oplus\cdots\oplus (q_{b}\circ f_{\ell})(x_{r_{\ell}})\big),
\end{multline*}
where 
\begin{compactitem}[\ \ $\circ$]
\item $q\colon\R^d\longrightarrow b$ is the orthogonal projection onto the $(d-m)$-dimensional subspace $b\in B$ of $\R^d$,
\item $\dist(\, \cdot\, ,\sk_{k_{\ell}}(\Delta_{N_{\ell}}))$ denotes the distance function to the $k_{\ell}$-skeleton of the simplex $\Delta_{N_{\ell}}$, and
\item $a(x_1,\ldots,x_{r_{\ell}})=\tfrac1{r_{\ell}}\big( \dist(x_1, \sk_{k_{\ell}}(\Delta_{N_{\ell}}))+\cdots+\dist(x_{r_{\ell}}, \sk_{k_{\ell}}(\Delta_{N_{\ell}}))\big)$.
\end{compactitem}

\noindent
Next we consider the compact subsets
\[
S_{\ell}:=\Phi_{\ell}^{-1}(E(\xi_{\ell})^{G_{\ell}})
\qquad\text{and}\qquad
T_{\ell}:=\Phi_{\ell}(S_{\ell}) =\im(\Phi_{\ell})\cap E(\xi_{\ell})^{G_{\ell}}.
\]
where $T_{\ell}\subseteq E(\xi_{\ell})^{G_{\ell}}\cong E(\gamma)$.
The set $S_{\ell}$ contains of all points $(b, (x_1,\ldots,x_{r_{\ell}}))\in B\times K_{\ell}$ such that
\[
\dist(x_1, \sk_{k_{\ell}}(\Delta_{N_{\ell}}))=\cdots = \dist(x_{r_{\ell}}, \sk_{k_{\ell}}(\Delta_{N_{\ell}}))
\qquad\text{and}\qquad
(q_{b}\circ f_{\ell})(x_1) =\cdots= (q_{b}\circ f_{\ell})(x_{r_{\ell}}).
\]
Since $(x_1,\ldots,x_{r_{\ell}})\in K_{\ell}$, then there exist $r_{\ell}$ pairwise disjoint faces $\sigma_1^{\ell},\ldots,\sigma_{r_{\ell}}^{\ell}$ of $\Delta_{N-{\ell}}$ such that
\[
(x_1,\ldots,x_{r_{\ell}})\in \relint\sigma_1^{\ell}\times \cdots\times \relint\sigma_{r_{\ell}}^{\ell},
\]
and at least one of the faces $\sigma_1^{\ell},\ldots,\sigma_{r_{\ell}}^{\ell}$ belongs to the $k_{\ell}$-skeleton of the simplex $\Delta_{N_{\ell}}$ \cite[Lem.\,4.2\,(iii)]{Blagojevic2014}.
Indeed, if this would no be true all the faces $\sigma_1^{\ell},\ldots,\sigma_{r_{\ell}}^{\ell}$ would  be at least $(k_{\ell}+1)$-dimensional, implying the following contradiction
    \[
    N_{\ell}+1=|\Delta_{N_{\ell}}|\geq |\sigma_1^{\ell}|+\dots+|\sigma_{r_{\ell}}^{\ell}|\geq r_{\ell}(k_{\ell}+2)\geq r_{\ell}\big(\lceil \tfrac{r_{\ell}-1}{r_{\ell}}d\rceil+2\big)\geq (r_{\ell}-1)(d+2)+2=N_{\ell}+2.
    \]
Therefore, at least one of the faces $\sigma_1^{\ell},\ldots,\sigma_{r_{\ell}}^{\ell}$ lies in $\sk_{k_{\ell}}(\Delta_{N_{\ell}})$ and consequently
\[
\dist(x_1, \sk_{k_{\ell}}(\Delta_{N_{\ell}}))=\cdots = \dist(x_{r_{\ell}}, \sk_{k_{\ell}}(\Delta_{N_{\ell}}))=0,
\]
implying that all the faces $\sigma_1^{\ell},\ldots,\sigma_{r_{\ell}}^{\ell}$ lies in $\sk_{k_{\ell}}(\Delta_{N_{\ell}})$.

Thus, in order to conclude the proof of Theorem \ref{th : transversal generalized Van Kamplen Flores} we need to show that
\[
\emptyset\neq T_0\cap\cdots\cap T_m\subseteq E(\gamma),
\]
and for that we would like to use Lemma~\ref{lemma : intersection lemma}.

\subsubsection{~~}
\label{subsec : 01}

First, let $0\leq\ell\leq m$ and let $e_{\ell}=0$.
Then $r_{\ell}=1$, $N_{\ell}=0$, $K_{\ell}=\Delta_{N_{\ell}}$ is a point, $G_{\ell}$ is the trivial group,  and $S_{\ell}=B\times K_{\ell}$.
Consider the commutative diagram induced by the bundle map $\Phi_{\ell}\colon B\times K_{\ell} \longrightarrow E(\xi_{\ell})$ and the corresponding diagram in cohomology:
\[
\xymatrix{
B\times  K_{\ell} = S_{\ell}\ar[rr]^{\Phi_{\ell}|_{S_{\ell}}} \ar[dr]_{p_1} && T_{\ell}\ar[dl]^{q_{\ell}|_{T_{\ell}}} &   H^*(B;\F_p)\cong H^*(S_{\ell};\F_p)  && H^*(T_{\ell};\F_p)\ar[ll]_-{(\Phi_{\ell}|_{S_{\ell}})^*}\\
& B & & &  H^*(B;\F_p).\ar[ul]_{p_1^*}\ar[ur]^{(q_{\ell}|_{T_{\ell}})^*}
}
\]
Since $K_{\ell}$ is a point the map $p_1^*$ induced by the projection $p_1$ is the identity map.
Consequently, the map
\[
(q_{\ell}|_{T_{\ell}})^*\colon H^*(B;\F_p)\longrightarrow  H^*(T_{\ell};\F_p).
\]
is an injection.

\subsubsection{~~}

Next, let $0\leq\ell\leq m$, and let $e_{\ell}>0$.
Now we apply Theorem \ref{th : BU} to the $G_{\ell}$-equivariant bundle map $\Phi_{\ell}\colon B\times K_{\ell} \longrightarrow E(\xi_{\ell})$.
In order to do so we check the necessary assumptions.
Since
\begin{compactitem}[\ \ $\circ$]
\item $G_{\ell}=(\Z_p)^{e_{\ell}}$ is an elementary abelian group,
\item $B=G_{d-m}(\R^d)$ is a connected space with the trivial $G_{\ell}$-action,
\item $q_{\ell}\colon E(\xi_{\ell}) \longrightarrow B$ is a $G_{\ell}$-equivariant vector bundle where all fibers carry the same $G_{\ell}$-representation,
\item $q_{\ell}|_{ E(\xi_{\ell})^{G_{\ell}}}\colon E(\xi_{\ell})^{G_{\ell}}\longrightarrow B$ is  the fixed-point subbundle with the $G_{\ell}$-invariant orthogonal complement subbundle $q_{\ell}|_{C_{\ell}}\colon C_{\ell}\longrightarrow B$, $(E(\xi_{\ell})=C_{\ell}\oplus E(\xi_{\ell})^{G_{\ell}})$,
\item $F_{\ell}$ is the fiber of the vector bundle $q_{\ell}|_{C_{\ell}}\colon C_{\ell}\longrightarrow B$ over the point $b\in B$,
\item $\pi_1(B)$ acts trivially on the cohomology of the sphere $H^*(S(F_{\ell});\F_p)$,
\item the Euler class $0\neq \alpha_{\ell}\in H^{(r_{\ell}-1)(d-m+1)}(G_{\ell};\F_p)$ of the vector bundle $F_{\ell}\longrightarrow \E G_{\ell}\times_{G_{\ell}} F_{\ell}\longrightarrow \B G_{\ell}$ does not vanish, more precisely 
		\[
		 \alpha_{\ell}= \Big( \prod_{(a_1,\ldots,a_{e_{\ell}})\in\F_p^{e_{\ell}}{\setminus}\{0\} }(a_1 t_1+\cdots +a_{e_{\ell}} t_{e_{\ell}})\Big)^{\tfrac{d-m+1}{2}}   ,
		\] 
\item $\ind_{G_{\ell}}^{\pt}(K_{\ell};\F_p)\subseteq H^{\geq (r_{\ell}-1)(d-m+1)+1}(G_{\ell};\F_p)$ because $K_{\ell}$ is $((r_{\ell}-1)(d-m+1)-1)$-connected,
\end{compactitem} 
we have that $\alpha_{\ell}\notin \ind_{G_{\ell}}^{\pt}(K_{\ell};\F_p)$ and Theorem \ref{th : BU} can be applied on the $G_{\ell}$-equivariant bundle map $\Phi_{\ell}\colon B\times K_{\ell} \longrightarrow E(\xi_{\ell})$. 
Thus, the following map in \v{C}ech cohomology induced by $q_{\ell}$ is injective:
\[
(q_{\ell}|_{T_{\ell}})^*\colon H^*(B;\F_p)\longrightarrow  H^*(T_{\ell};\F_p).
\]

\subsubsection{~~}
Finally, Lemma~\ref{lemma : intersection lemma} comes into play.
Since,
\begin{compactitem}[\ \ $\circ$]
\item $T_{\ell}$ is a compact subset of 	$E(\gamma)$ for every $0\leq\ell\leq m$,
\item $(q_{\ell}|_{T_{\ell}})^*\colon H^*(B;\F_p)\longrightarrow  H^*(T_{\ell};\F_p)$ is injective for every $0\leq\ell\leq m$, and
\item $0\neq e(\gamma)^m \in H^{(d-m)m}(B;\F_p)$ according to $p(d-m)$ being even and Lemma \ref{lemma : euler class},
\end{compactitem}
we can apply Lemma~\ref{lemma : intersection lemma} and get that
\[
T_0\cap\cdots\cap T_m\neq\emptyset.
\]
This concludes the proof of Theorem~\ref{th : transversal generalized Van Kamplen Flores}.\hfill\qed

\subsection{Proof of the topological transversal weak colored Tverberg theorem}
\label{subsec : WCT}
Let $B:=G_{d-m}(\R^d)$ be the Grassmann manifold, and let $\gamma:=\gamma^{d-m}(\R^d)$ be the canonical bundle.
Without lost of generality we can assume that $m\geq 1$.
The proof of Theorem \ref{th : topological transversal weak colored tverberg theorem} is done in the footsteps of the proof of Theorem~\ref{th : transversal generalized Van Kamplen Flores}.
The only difference will occur in the definition of the bundles $\tau_{\ell}$ and consequently bundle maps $\Phi_{\ell}$.

\subsubsection{~~}

Again, fix an integer $0\leq\ell\leq m$, and define $K_{\ell}:=(\Delta_{N_{\ell}})^{\times r_{\ell}}_{\Delta(2)}$. 
As we have seen the complex $K_{\ell}$ is an $(N_{\ell}-r_{\ell}+1)$-dimensional and $(N_{\ell}-r_{\ell})$-connected CW complex.
The symmetric group $\Sym_{r_{\ell}}$ acts freely on $K_{\ell}$ by permuting factors in the product.

The regular embedding $\mathrm{reg}\colon (\Z/p)^{e_{\ell}}\longrightarrow \Sym_{r_{\ell}}$ of the elementary abelian group $(\Z/p)^{e_{\ell}}$
identifies elementary abelian group $G_{\ell}:=(\Z/p)^{e_{\ell}}$ with a subgroup $\im(\mathrm{reg})$ of the symmetric group $\Sym_{r_{\ell}}$.

Once more, $\R^{r_{\ell}}$ is a vectors space with the (left) action of the symmetric group $\Sym_{r_{\ell}}$ given by permutation of coordinates. 
The subspace $W_{r_{\ell}}:=\{(t_1,\ldots,t_{r_{\ell}})\in \R^{r_{\ell}} : \sum t_i=0\}$ is a $\Sym_{r_{\ell}}$-invariant subspace, and $G_{\ell}$ acts on both $\R^{r_{\ell}}$ and $W_{r_{\ell}}$ via the regular embedding.

Let $\tau_{\ell}$ be the the trivial vector bundle $B\times W_{r_{\ell}}^{\oplus d+1}\longrightarrow B$.
The action of $G_{\ell}$ on $W_{r_{\ell}}^{\oplus d+1}$ is diagonal and makes $\tau_{\ell}$ into a $G_{\ell}$-equivariant vector bundle.
As we have seen,  $\gamma^{\oplus r_{\ell}}$ is also a $G_{\ell}$-equivariant vector bundle.
Thus the vector bundle 
\begin{equation}
	\label{def of xi - 02}
	\xi_{\ell}:=\tau_{\ell}\oplus \gamma^{\oplus r_{\ell}}
\end{equation}
inherits the structure of a $G_{\ell}$-equivariant vector bundle via the diagonal action.
Since $(W_{r_{\ell}}^{\oplus d+1})^{G_{\ell}}=\{0\}$, the fixed point set of the total space of $\xi_{\ell}$ is
\[
 E(\xi_{\ell})^{G_{\ell}}=
 E(\tau_{\ell}\oplus \gamma^{\oplus r_{\ell}})^{G_{\ell}}\cong 
 E(\gamma^{\oplus r_{\ell}})^{G_{\ell}}\cong
 E(\gamma).
\]
 
\subsubsection{~~}
\label{subsubsec: definition of function WCT}

The vertices of the simplex $\Delta_{N_{\ell}}$ are colored by $d+1$ colors, where each color class has cardinality at most $2r_{\ell}-1$.
Set $\vertices(\Delta_N)=C_0\sqcup\cdots\sqcup C_{d}$ where $|C_i|\leq 2r_{\ell}-1$ for all $0\leq i\leq d$.
Following the idea from \cite[Lem.\,4.2\,(ii)]{Blagojevic2014} we define $\Sigma_i^{\ell}$, $0\leq i\leq d$, to be the subcomplex of $\Delta_{N_{\ell}}$ consisting of all faces with at most one vertex in $C_i$.
Then the rainbow subcomplex $R_{N_{\ell}}$ coincides with the intersection $\Sigma_0^{\ell}\cap\cdots\cap\Sigma_d^{\ell}$.

Now we define a continuous $G_{\ell}$-equivariant bundle map $\Phi_{\ell}\colon B\times K_{\ell} \longrightarrow E(\xi_{\ell})$ as follows: For the point $(b,(x_1,\ldots,x_{r_{\ell}}))\in B\times K_{\ell}$ let
\begin{eqnarray*}
\Phi_{\ell} (b,(x_1,\ldots,x_{r_{\ell}})) &:=& 	 \big(b,\dist(x_1,\Sigma_0^{\ell})-a_0(x_1,\ldots,x_{r_{\ell}}),\ldots,  \dist(x_{r_{\ell}},\Sigma_0^{\ell})-a_0(x_1,\ldots,x_{r_{\ell}}) \big) \oplus\\
& & \qquad\qquad\qquad\qquad\qquad\quad\qquad \cdots \\
& &  \big(b,\dist(x_1,\Sigma_d^{\ell})-a_d(x_1,\ldots,x_{r_{\ell}}),\ldots,  \dist(x_{r_{\ell}},\Sigma_d^{\ell})-a_d(x_1,\ldots,x_{r_{\ell}}) \big) \oplus\\
& & \big((q_{b}\circ f_{\ell})(x_1) \oplus\cdots\oplus (q_{b}\circ f_{\ell})(x_{r_{\ell}})\big)
\end{eqnarray*}
where 
\begin{compactitem}[\ \ $\circ$]
\item $q\colon\R^d\longrightarrow b$ is the orthogonal projection onto the $(d-m)$-dimensional subspace $b\in B$ of $\R^d$,
\item $\dist(\, \cdot\, , \Sigma_i^{\ell})$ denotes the distance function to the subcomplex $\Sigma_i^{\ell}$ where $0\leq i\leq d$, and
\item $a_i(x_1,\ldots,x_{r_{\ell}})=\tfrac1{r_{\ell}}\big( \dist(x_1,\Sigma_i^{\ell})+\cdots+\dist(x_{r_{\ell}},\Sigma_i^{\ell}))$ for $0\leq i\leq d$.
\end{compactitem}

\noindent
Again, we consider the compact subsets
\[
S_{\ell}:=\Phi_{\ell}^{-1}(E(\xi_{\ell})^{G_{\ell}})
\qquad\text{and}\qquad
T_{\ell}:=\Phi_{\ell}(S_{\ell}) =\im(\Phi_{\ell})\cap E(\xi_{\ell})^{G_{\ell}}.
\]
where $T_{\ell}\subseteq E(\xi_{\ell})^{G_{\ell}}\cong E(\gamma)$.
The set $S_{\ell}$ contains of all points $(b, (x_1,\ldots,x_{r_{\ell}}))\in B\times K_{\ell}$ such that
\[
\dist(x_1,\Sigma_0^{\ell})=\cdots=\dist(x_{r_{\ell}},\Sigma_0^{\ell}), \quad \ldots \quad  ,
\dist(x_1,\Sigma_d^{\ell})=\cdots=\dist(x_{r_{\ell}},\Sigma_d^{\ell}),
\]
and
\[
(q_{b}\circ f_{\ell})(x_1) =\cdots= (q_{b}\circ f_{\ell})(x_{r_{\ell}}).
\]
Since the point $(x_1,\ldots,x_{r_{\ell}})\in K_{\ell}$, then we can find $r_{\ell}$ unique pairwise disjoint faces $\sigma_1^{\ell},\ldots,\sigma_{r_{\ell}}^{\ell}$ with the property that
\[
(x_1,\ldots,x_{r_{\ell}})\in \relint\sigma_1^{\ell}\times \cdots\times \relint\sigma_{r_{\ell}}^{\ell}.
\]
Moreover, for every $i$ in the range $0\leq i\leq d$ there exists at least one of the faces $\sigma_1^{\ell},\ldots,\sigma_{r_{\ell}}^{\ell}$ that belongs to the subcomplex $\Sigma_i$ of the simplex $\Delta_{N_{\ell}}$, \cite[Lem.\,4.2\,(ii)]{Blagojevic2014}.
Indeed, if this would no be true for an index $i$ then all the faces $\sigma_1^{\ell},\ldots,\sigma_{r_{\ell}}^{\ell}$ would have at least two vertices from the set $C_i$ and we have the contradiction
	\[
		2r_{\ell}-1\geq |C_i|\geq |\sigma_1^{\ell}\cap C_i|+\cdots + |\sigma_{r_{\ell}}^{\ell}\cap C_i|\geq 2r_{\ell}.
	\]
Thus for every index $i$ at least one of the faces $\sigma_1^{\ell},\ldots,\sigma_{r_{\ell}}^{\ell}$ lies in $\Sigma_i$ and consequently
\[
\dist(x_1,\Sigma_0^{\ell})=\cdots=\dist(x_{r_{\ell}},\Sigma_0^{\ell})=0, \quad \ldots \quad  ,
\dist(x_1,\Sigma_d^{\ell})=\cdots=\dist(x_{r_{\ell}},\Sigma_d^{\ell})=0,
\]
implying that all the faces $\sigma_1^{\ell},\ldots,\sigma_{r_{\ell}}^{\ell}$ lie in every subcomplex $\Sigma_i$, meaning they belong to the intersection $\Sigma_0^{\ell}\cap\cdots\cap\Sigma_d^{\ell}$ -- the rainbow subcomplex. 

Therefore, in order to finalize the proof of Theorem \ref{th : topological transversal weak colored tverberg theorem} we need to show, as in the proof of the previous theorem, that
\[
\emptyset\neq T_0\cap\cdots\cap T_m\subseteq E(\gamma).
\]

\subsubsection{~~}

Again first assume that $0\leq\ell\leq m$ and that $e_{\ell}=0$.
We proceed as in Section~\ref{subsec : 01}. 
Now $r_{\ell}=1$, $N_{\ell}=0$, $K_{\ell}=\Delta_{N_{\ell}}$ is a point, $G_{\ell}$ is the trivial group,  and $S_{\ell}=B\times K_{\ell}$.
We consider the commutative diagrams induced by the bundle map $\Phi_{\ell}\colon B\times K_{\ell} \longrightarrow E(\xi_{\ell})$:
\[
\xymatrix{
B\times  K_{\ell} = S_{\ell}\ar[rr]^{\Phi_{\ell}|_{S_{\ell}}} \ar[dr]_{p_1} && T_{\ell}\ar[dl]^{q_{\ell}|_{T_{\ell}}} &   H^*(B;\F_p)\cong H^*(S_{\ell};\F_p)  && H^*(T_{\ell};\F_p)\ar[ll]_-{(\Phi_{\ell}|_{S_{\ell}})^*}\\
& B & & &  H^*(B;\F_p).\ar[ul]_{p_1^*}\ar[ur]^{(q_{\ell}|_{T_{\ell}})^*}
}
\]
The map $p_1^*$ induced by the projection $p_1$ is the identity map.
Hence, the map in cohomology 
\[
(q_{\ell}|_{T_{\ell}})^*\colon H^*(B;\F_p)\longrightarrow  H^*(T_{\ell};\F_p).
\]
is an injection.

\subsubsection{~~}

Let $0\leq\ell\leq m$ and $e_{\ell}>0$. 
We apply Theorem \ref{th : BU} to the $G_{\ell}$-equivariant bundle map $\Phi_{\ell}\colon B\times K_{\ell} \longrightarrow E(\xi_{\ell})$.
For that we check the necessary assumptions.
Since
\begin{compactitem}[\ \ $\circ$]
\item $G_{\ell}=(\Z_p)^{e_{\ell}}$ is an elementary abelian group,
\item $B=G_{d-m}(\R^d)$ is a connected space with the trivial $G_{\ell}$-action,
\item $q_{\ell}\colon E(\xi_{\ell}) \longrightarrow B$ is a $G_{\ell}$-equivariant vector bundle where all fibers carry the same $G_{\ell}$-representation,
\item $q_{\ell}|_{ E(\xi_{\ell})^{G_{\ell}}}\colon E(\xi_{\ell})^{G_{\ell}}\longrightarrow B$ is  the fixed-point subbundle with the $G_{\ell}$-invariant orthogonal complement subbundle $q_{\ell}|_{C_{\ell}}\colon C_{\ell}\longrightarrow B$, $(E(\xi_{\ell})=C_{\ell}\oplus E(\xi_{\ell})^{G_{\ell}})$,
\item $F_{\ell}$ is the fiber of the vector bundle $q_{\ell}|_{C_{\ell}}\colon C_{\ell}\longrightarrow B$ over the point $b\in B$,
\item $\pi_1(B)$ acts trivially on the cohomology of the sphere $H^*(S(F_{\ell});\F_p)$,
\item the Euler class $0\neq \alpha_{\ell}\in H^{(r_{\ell}-1)(2d-m+1)}(G_{\ell};\F_p)$ of the vector bundle $F_{\ell}\longrightarrow \E G_{\ell}\times_{G_{\ell}} F_{\ell}\longrightarrow \B G_{\ell}$ does not vanish, more precisely 
		\[
		 \alpha_{\ell}= \Big( \prod_{(a_1,\ldots,a_{e_{\ell}})\in\F_p^{e_{\ell}}{\setminus}\{0\} }(a_1 t_1+\cdots +a_{e_{\ell}} t_{e_{\ell}})\Big)^{\tfrac{2d-m+1}{2}},
		\] 
\item $\ind_{G_{\ell}}^{\pt}(K_{\ell};\F_p)\subseteq H^{\geq (r_{\ell}-1)(2d-m+1)+1}(G_{\ell};\F_p)$ because $K_{\ell}$ is $((r_{\ell}-1)(2d-m+1)-1)$-connected,
\end{compactitem} 
we conclude that $\alpha_{\ell}\notin \ind_{G_{\ell}}^{\pt}(K_{\ell};\F_p)$, and therefore Theorem \ref{th : BU} can be applied on the $G_{\ell}$-equivariant bundle map $\Phi_{\ell}\colon B\times K_{\ell} \longrightarrow E(\xi_{\ell})$. 
Consequently, the following map induced by $q_{\ell}$ is injective:
\[
(q_{\ell}|_{T_{\ell}})^*\colon H^*(B;\F_p)\longrightarrow  H^*(T_{\ell};\F_p).
\]

\subsubsection{~~}
In the final step we apply Lemma~\ref{lemma : intersection lemma}.
Since,
\begin{compactitem}[\ \ $\circ$]
\item $T_{\ell}$ is a compact subset of 	$E(\gamma)$ for every $0\leq\ell\leq m$,
\item $(q_{\ell}|_{T_{\ell}})^*\colon H^*(B;\F_p)\longrightarrow  H^*(T_{\ell};\F_p)$ is injective for every $0\leq\ell\leq m$, and
\item $0\neq e(\gamma)^m \in H^{(d-m)m}(B;\F_p)$ does not vanish according to $p(d-m)$ being even and Lemma \ref{lemma : euler class},
\end{compactitem}
we can apply Lemma~\ref{lemma : intersection lemma} and get that
\[
T_0\cap\cdots\cap T_m\neq\emptyset.
\]
This concludes the proof of Theorem~\ref{th : topological transversal weak colored tverberg theorem}. \hfill\qed

\providecommand{\noopsort}[1]{}
\providecommand{\bysame}{\leavevmode\hbox to3em{\hrulefill}\thinspace}
\providecommand{\MR}{\relax\ifhmode\unskip\space\fi MR }
\providecommand{\MRhref}[2]{%
  \href{http://www.ams.org/mathscinet-getitem?mr=#1}{#2}
}
\providecommand{\href}[2]{#2}

\end{document}